\documentclass[12pt,twoside,a4paper]{amsart}
\usepackage{amsfonts,latexsym,rawfonts,amsmath,amssymb,amsthm,mathrsfs,mathtools}
\usepackage[colorlinks=true,urlcolor=blue,linkcolor=blue,anchorcolor=blue,citecolor=blue]{hyperref}
\usepackage{graphicx}
\usepackage{comment}
\usepackage[]{caption2}

\usepackage[pagewise]{lineno}
\usepackage{lipsum}
\usepackage[raggedright]{sidecap}\sidecaptionvpos{figure}{t} 
\usepackage{tikz}
\usepackage{enumerate}

\allowdisplaybreaks[4]
\usepackage{hyperref}

\makeatletter
\@namedef{subjclassname@2020}{\textup{2020} Mathematics Subject Classification}
\makeatother

\providecommand{\bysame}{\leavevmode\hbox to3em{\hrulefill}\thinspace}
\providecommand{\MR}{\relax\ifhmode\unskip\space\fi MR }

\providecommand{\href}[2]{#2}

\theoremstyle{plain}
\newtheorem{theorem}{Theorem}[section]
\newtheorem{corollary}[theorem]{Corollary}
\newtheorem{lemma}[theorem]{Lemma}
\newtheorem{proposition}[theorem]{Proposition}

\theoremstyle{definition}
\newtheorem{definition}[theorem]{Definition}
\newtheorem{conditions}{Conditions}
\newtheorem{example}{Example}

\theoremstyle{remark}
\newtheorem*{remark*}{Remark}

\numberwithin{equation}{section}

\begin{document}

\title{Noncollapsing for Curvature Flows with Inhomogeneous Speeds}

\author{Weimin Sheng}
\address{Weimin Sheng: School of Mathematical Sciences, Zhejiang University, Hangzhou 310058, China.}
\email{weimins@zju.edu.cn}
\author{Ye Zhu}
\address{Ye Zhu: School of Mathematical Sciences, Zhejiang University, Hangzhou 310058, China.}
\email{ye.zhu@zju.edu.cn}


\maketitle

\begin{abstract}
	We study closed, embedded hypersurfaces in Euclidean space evolving by fully nonlinear curvature flows, whose speed is given by a symmetric, monotone increasing, $1$-homogeneous, positive underlying speed function $F$ composed with a modulating function $\Psi$. Under the assumption that $F$ is convex or inverse-concave and that $\Psi$ satisfies the corresponding structural conditions, we establish exterior noncollapsing estimates for the flow. The main difficulty stems from the nonlinearity of the evolution equation satisfied in the viscosity sense by the exscribed curvature, whereas in previous works it is a solution to the linearized flow. Moreover, in the case where $F$ is inverse-concave, we adopt an argument different from that of Andrews and Langford to treat off-diagonal extrema.
\end{abstract}

\setlength{\baselineskip}{16pt}
\setlength{\parskip}{3pt}


\section{Introduction}

We consider a smooth one-parameter family of embeddings 
\[
X(\cdot, t): M^n \to \mathbb{R}^{n+1}, \quad t\in [0,T), 
\]
where $M^n$ is a closed $n$-dimensional manifold, $n\geq 2$, and $\mathbb{R}^{n+1}$ is equipped with the standard Euclidean metric. 
The image $M_t \coloneq X(M,t)$ is a smooth, oriented hypersurface with unit normal vector field $\nu(\cdot,t)$.
We denote by 
\[
\kappa(\cdot,t) = (\kappa_1(\cdot,t), \ldots, \kappa_n(\cdot,t))
\]
the $n$-tuple of principal curvatures of $M_t$ with respect to $\nu(\cdot,t)$.

The evolution of $M_t$ is governed by the fully nonlinear curvature flow
\begin{equation*} 
	\partial_t X(x,t) = -\Psi(F(x,t))\, \nu(x,t),
	\tag{CF} \label{CF}
\end{equation*}
where 
\begin{itemize}
	\item The underlying speed function $F$ is given by 
	\[
	F(x,t) = f(\kappa(x,t)), 
	\]
	for some smooth function $f: \Gamma \subset \mathbb{R}^n \to \mathbb{R}$, defined on an open, symmetric cone $\Gamma$ containing the positive cone $\Gamma_+$. 
	Throughout our main results, we assume that $f$ satisfies the following structural properties: 
	\begin{conditions}\label{conditions_f}
		\begin{itemize}
			\item[\rm{(i)}] \textbf{Symmetry:} $f$ is invariant under permutation of its variables. 
			\item[\rm{(ii)}] \textbf{Parabolicity:} $f$ is monotone increasing in each variable.
			\item[\rm{(iii)}] \textbf{Homogeneity:} $f$ is homogeneous of degree one.
			\item[\rm{(iv)}] \textbf{Positivity:} $f>0$ on $\Gamma$.
		\end{itemize}
	\end{conditions}
	\noindent Whenever $f$ satisfies Conditions~\ref{conditions_f}(i)--(iv), we refer to it as an \textit{admissible underlying speed function} for \eqref{CF}. Typical examples of admissible underlying speed functions include the power means $\left(\frac{1}{n} \sum_{i=1}^n \lambda_i^r \right)^{\frac{1}{r}}$, with $r\neq 0$, defined on the positive cone $\Gamma_+$; the roots of the elementary symmetric polynomials $\sigma_k^{\frac{1}{k}}$, for $1 \leq k\leq n$, defined on the cone $\Gamma_k\coloneq \{\lambda\in \mathbb{R}^{n}: \sigma_j(\lambda)>0, \, \forall \, 1\leq j\leq k\}$; and the roots of the ratios of elementary symmetric polynomials $(\sigma_k/ \sigma_l)^{\frac{1}{k-l}}$, for $0\leq l<k \leq n$, also defined on $\Gamma_k$. 
	\item The modulating function $\Psi: (0, +\infty) \to \mathbb{R}$ is assumed to be $C^2$, and to satisfy certain (not necessarily all) of the structural conditions listed below, depending on the context of application: 
	\begin{conditions}\label{conditions_psi}
		\begin{itemize}
			\item[\rm{(i)}] $\Psi'(s)>0, \ \forall s>0$. 
			\item[\rm{(ii)}] Either $\mathrm{(a)} \ \Psi'(s)s-\Psi(s) \leq 0, \ \forall s>0$, or $\mathrm{(b)} \ \Psi'(s)s-\Psi(s) \geq 0, \ \forall s>0$.
			\item[\rm{(iii)}] Either $\mathrm{(a)} \ \Psi''(s) \geq 0, \ \forall s>0$, or $\mathrm{(b)} \ \Psi''(s) \leq 0, \ \forall s>0$. 
			\item[\rm{(iv)}] $\Psi''(s)s+2\Psi'(s) \geq 0, \ \forall s>0$.
		\end{itemize}
	\end{conditions}
\end{itemize}

\bigskip

In addition to the extensive literature on homogeneous curvature flows, there has been growing interest in flows driven by inhomogeneous speeds. Chow and Tsai~\cite{chow1997expansion} studied expanding convex hypersurfaces whose speeds are given by inhomogeneous functions of the principal radii, establishing long-time existence and asymptotic roundness after rescaling. In the contraction setting, Li and Lv~\cite{li2020contracting} and McCoy~\cite{mccoy2021contraction} considered flows in space forms with general inhomogeneous speeds, proving that, under suitable structural conditions, the hypersurfaces shrink to a point in finite time, and that the rescaled hypersurfaces converge smoothly to the round sphere. These results indicate that the flow~\eqref{CF} shares essential structural features with classical homogeneous flows, while also presenting distinctive analytic challenges. 

The analysis of curvature flows has been profoundly influenced by the development of noncollapsing estimates, which provide a quantitative measure of embeddedness and play a crucial role in understanding singularity formation. The concept originates from joint work of the first author and Wang~\cite{sheng2009singularity}, where they introduced the notion of $\delta$-noncollapsing for the mean-convex mean curvature flow and proved that each point of the evolving hypersurface admits an interior touching ball of radius $\delta/H$, for some constant $\delta>0$. Andrews~\cite{andrews2012noncollapsing} later gave a remarkably simple and flexible proof based on a two-point maximum principle. In this approach, one considers the function
\[k(x,y,t)\coloneq \frac{2\left< X(x,t)-X(y,t), \nu(x,t) \right>}{|X(x,t)-X(y,t)|^2}, \]
and defines the inscribed and exscribed curvatures
\[\overline k(x,t)\coloneq \sup_{y\in M\setminus\{x\}}k(x,y,t), \quad \underline k(x,t)\coloneq \inf_{y\in M\setminus\{x\}}k(x,y,t); \]
see Section~\ref{subsec:ballcurvatures} for a detailed exposition and further geometric properties. Andrews showed that the ratio $\overline k/H$ is monotone non-increasing along the mean curvature flow, thereby providing a conceptually transparent proof of noncollapsing. Building on this idea, Brendle~\cite{brendle2015sharp} established sharp estimates via Stampacchia iteration. More precisely, given any $\varepsilon>0$, the inscribed curvature $\overline k$ is bounded from above by $(1+\varepsilon)H$ at each point where the curvature is sufficiently large. Haslhofer and Kleiner~\cite{haslhofer2015brendle} subsequently found an alternative proof by a blow-up argument, while Langford~\cite{langford2015optimal} improved the result in the setting of $(m+1)$-convex mean curvature flow for $m=0, \ldots, n-2$. 

Parallel progress has been made for fully nonlinear curvature flows. Andrews, Langford, and McCoy~\cite{andrews2013non} first established noncollapsing estimates for curvature flows with concave or convex speeds. Later, Andrews et al.~\cite{andrews2015non} extended these results to space forms with a suitable modification. A further breakthrough was achieved by Andrews and Langford~\cite{andrews2016two}, who introduced a two-sided noncollapsing condition and developed a maximum principle framework that applies to a wide class of fully nonlinear curvature flows, thereby unifying and generalizing previous approaches. More recently, Langford and Lynch~\cite{langford2020sharp} refined the analysis using Stampacchia iteration to derive sharp one-sided estimates --- including cylindrical estimates, inscribed curvature pinching for concave speeds, and exscribed curvature pinching for convex speeds --- leading to classification results for ancient solutions and rigidity theorems characterizing shrinking spheres. Collectively, these contributions form a coherent theory that both parallels and extends the picture established in the mean curvature case. 

The impact of the two-point maximum principle reaches far beyond analytic estimates. This technique was a key ingredient in Brendle's resolution of the Lawson conjecture~\cite{brendle2013embedded} and in the Andrews--Li's classification of embedded constant mean curvature tori in the three-sphere~\cite{andrews2015embedded}, which confirmed the Pinkall--Sterling conjecture. Moreover, the sharp noncollapsing estimate was applied by Brendle and Huisken~\cite{brendle2016mean} in the construction of mean curvature flow with surgery for mean-convex surfaces in $\mathbb{R}^3$, effectively replacing the cylindrical estimate employed in the earlier work of Huisken and Sinestrari~\cite{huisken2009mean}. In subsequent work, suitable noncollapsing estimates were established to derive the pointwise curvature derivative estimate, enabling the extension of the surgery procedure to fully nonlinear curvature flows for two-convex hypersurfaces in Riemannian manifolds~\cite{brendle2017fully}. These examples illustrate how noncollapsing estimates not only sharpen analytic control of curvature flows but also open the door to resolving long-standing geometric problems. 

Within this broader landscape, it is natural to ask whether analogous estimates can be developed for more general flows (relaxing the restriction of $1$-homogeneity). Addressing this question provides the starting point of the present work.

\bigskip

We now state our main results concerning the exscribed curvature $\underline k$. It is worth noting that the conditions imposed below are sufficient for our results, though alternative assumptions may also yield similar conclusions. 

\begin{theorem}\label{thm_convex}
	Suppose that $f: \Gamma \to \mathbb{R}$ is an admissible underlying speed function and $\Psi$ satisfies \textnormal{Conditions~\ref{conditions_psi}(i)}. Let $X: M^n \times [0,T) \to \mathbb{R}^{n+1}$ be an embedded solution of \textnormal{\eqref{CF}}. If $f$ is convex and $\Psi$ satisfies \textnormal{Conditions~\ref{conditions_psi}(iia)(iiia)}, then $X$ is exterior noncollapsing; that is, for all $(x,t)\in M \times [0,T)$,
	\[\underline k(x,t)\geq -\beta F(x,t), \]
	where $\beta>0$ is a constant such that $\underline k(x,0)\geq -\beta F(x,0)$ for all $x\in M$. 
\end{theorem}

\begin{example}\label{example_psi_1}
	Apart from the trivial case $\Psi(s)=s$, examples satisfying \textnormal{Conditions~\ref{conditions_psi}(i)(iia)(iiia)} include: \\
	\textup{(i)} $\Psi(s)=\sqrt{s^2+1}$, \quad
	\textup{(ii)} $\Psi(s)=\log(1+e^s)$, \quad 
	\textup{(iii)} $\Psi(s)=\log(e^s+e^{-s})$. \\
	In particular, no function satisfying these conditions attains negative values. 
\end{example}

\begin{theorem}\label{thm_inverse-concave}
	Suppose that $f: \Gamma_+ \to \mathbb{R}$ is an admissible underlying speed function and $\Psi$ satisfies \textnormal{Conditions~\ref{conditions_psi}(i)}. Let $X: M^n \times [0,T) \to \mathbb{R}^{n+1}$ be an embedded solution of \textnormal{\eqref{CF}}. If $f$ is inverse-concave (see Definition~\ref{def_inverse-concave} below for the precise definition) and $\Psi$ satisfies \textnormal{Conditions~\ref{conditions_psi}(iib)(iiib)(iv)}, then for all $(x,t)\in M \times [0,T)$, we have
	\[\frac{\underline k(x,t)}{F(x,t)}\geq c_0 \]
	with $c_0=\min_{M \times \{0\}}\frac{\underline k}{F}>0$.
\end{theorem}

\begin{example}\label{example_psi_2}
	Representative examples of $\Psi$ satisfying \textnormal{Conditions~\ref{conditions_psi}(i)(iib)(iiib)(iv)} include: \\
	\textup{(i)} $\Psi(s)=-s^{-\alpha}, \ 0<\alpha\leq 1$, \quad
	\textup{(ii)} $\Psi(s)=-\log(1+1/s)$, \quad
	\textup{(iii)} $\Psi(s)=-\log(1+s)/s$, \\
	\textup{(iv)} $\Psi(s)=-\arctan(1/s)$, \quad
	\textup{(v)} $\Psi(s)=s-e^{-s}$ \ (sign-changing). \\
	Note that the only positive function satisfying these conditions is $\Psi(s)=s$ (up to normalization). 
\end{example}

\bigskip

We briefly recall the definitions of the inradius $\rho_-$ and the circumradius $\rho_+$. 
Let $\Omega_t$ denote the region enclosed by $M_t$. Then  
\[
\begin{aligned}
	\rho_-(t) & \coloneq \sup\{ r>0 : B_r(y) \subset \Omega_t \text{ for some } y\in \mathbb{R}^{n+1} \}, \\
	\rho_+(t) & \coloneq \inf\{ r>0 : \Omega_t \subset B_r(y) \text{ for some } y\in \mathbb{R}^{n+1}\}.
\end{aligned}
\]
As a direct consequence of the exterior noncollapsing estimate, we obtain a uniform bound on the ratio of the circumradius to the inradius. 
\begin{corollary}\label{cor}
	Under the same assumptions as in Theorem~\ref{thm_inverse-concave}, suppose further that either $n=2$, or 
	\[\frac{\min\{\lambda_1, \dots, \lambda_n\}}{f(\lambda_1, \dots, \lambda_n)} \text{ approaches zero on the boundary of } \Gamma_+. \]
	Then, for all $t\in [0,T)$, the following estimate holds: 
	\[\rho_+(t) \leq K_0 \rho_-(t), \]
	where $K_0>0$ depends only on the underlying speed function $F$ and the initial hypersurface $M_0$. 
\end{corollary}

\bigskip

Before outlining the structure of the paper, we note that, due to technical difficulties, no interior noncollapsing estimate is currently known for the flow~\eqref{CF} with a nonlinear function~$\Psi$; the only established case is the trivial one~\cite[Corollary~3]{andrews2013non}.

\bigskip

The paper is organized as follows. In Section~\ref{sec:pre}, we recall the characterization of inverse-concavity and review the geometric properties of the ball curvatures. Section~\ref{sec:aux} is devoted to several auxiliary computations that prepare for the analysis of the ``interior case" in the subsequent proofs. In Section~\ref{sec:evo_underlinek}, we establish a key intermediate step toward the main theorems: the viscosity supersolution satisfied by the exscribed curvature $\underline k$. In particular, we provide an alternative proof of Lemma~\ref{inverse-concave_interior}, originally due to Andrews and Langford~\cite[Proposition~2.3]{andrews2016two}. Finally, the proofs of our main results are presented in Sections~\ref{sec:proof_convex} and \ref{sec:proof_inverse-concave}.

\bigskip

\section{Preliminaries}\label{sec:pre}

\subsection{Notation and evolution equations}

We begin by fixing notation and recalling the evolution equations for the geometric quantities associated with the flow~\eqref{CF}, which will be used extensively in the sequel. Let $g$ denote the induced metric on $M^n$, and let $\nabla$ be the Levi-Civita connection. We write $h_{ij}$ for the components of the second fundamental form of $M$.

Throughout this paper, $\dot F^{ij}$ denotes the derivative of $F$ with respect to $h_{ij}$, defined by 
\[
\dot F^{ij}(A)B_{ij} \coloneq \frac{d}{ds}\bigg|_{s=0}F(A+sB)
\]
for symmetric matrices $A$ and $B$. 
For brevity, we introduce the following notations:  
\[
\Delta_F \coloneq \dot F^{ij} \nabla_i\nabla_j, 
\quad |A|_F^2 \coloneq \dot F^{ij} {h^2}_{ij} = \dot F^{ij} g^{kl}h_{ik} h_{jl}, 
\quad |\nabla F|_F^2 \coloneq \dot F^{ij}\nabla_i F\nabla_j F.
\]

By a direct computation, one obtains the following evolution equations under the flow~\eqref{CF}:
\begin{itemize}
	\item[(i)] $\partial_t g_{ij}= -2\Psi h_{ij}$.
	\item[(ii)] $\partial_t \nu= \nabla \Psi$.
	\item[(iii)] $\partial_t h_{ij}=\Psi' \Delta_F h_{ij}+ (\Psi' \ddot F^{pq,rs}+ \Psi'' \dot F^{pq} \dot F^{rs}) \nabla_i h_{pq} \nabla_j h_{rs} + \Psi' |A|_F^2 h_{ij} - (\Psi' F + \Psi){h^2}_{ij}$.
\end{itemize}

Moreover, the underlying speed function $F$ satisfies
\begin{equation}\label{evo_F}
	\partial_t F = \Psi' \Delta_F F + \Psi''|\nabla F|_F^2 + |A|_F^2 \Psi.
\end{equation}

\bigskip

\subsection{Inverse-concave functions}

In this subsection, we introduce several auxiliary results that will be used later. 
While our main theorems require that $f$ satisfies all the structural properties of Conditions~\ref{conditions_f}, some of the lemmas below hold under weaker assumptions, where only certain properties are required.

\medskip

We begin by recalling the notion of inverse-concavity, which will play a crucial role in our analysis. 

\begin{definition}[Inverse-concavity]\label{def_inverse-concave}
	Let $\Gamma_+ \subset \mathbb{R}^n$ be the positive cone and let $f: \Gamma_+ \to \mathbb{R}$ be a function satisfying \textnormal{Conditions~\ref{conditions_f}(i)(iv)}. 
	We say that $f$ is inverse-concave if the dual function $f_*: \Gamma_+ \to \mathbb{R}$ defined by 
	\[
	f_*(\lambda_1^{-1}, \ldots, \lambda_n^{-1})\coloneq f(\lambda_1, \ldots, \lambda_n)^{-1}
	\]
	is concave.
\end{definition}

We refer to~\cite{andrews2013contracting} for a variety of explicit examples of inverse-concave functions.

\bigskip

For clarity, we fix notation for derivatives of $f$ with respect to its arguments $\lambda_i$: 
\[\dot f^i \coloneq \frac{\partial f}{\partial \lambda_i}, \quad \ddot f^{ij}\coloneq \frac{\partial^2 f}{\partial \lambda_i \partial \lambda_j}. \]

The following lemma provides a local characterization of inverse-concavity. 

\begin{lemma}[\cite{andrews2007pinching, andrews2013contracting}]\label{inverse-concave_origin}
	Let $f: \Gamma_+ \to \mathbb{R}$ be a function satisfying \textnormal{Conditions~\ref{conditions_f}(i)(iv)}. Then $f$ is inverse-concave if and only if 
	\[\ddot f^{ij}- \frac{2}{f} \dot f^i \dot f^j + 2\frac{\dot f^i}{\lambda_i}\delta_{ij} \]
	is nonnegative definite for all $\lambda \in \Gamma_+$ and, in addition, 
	\[\frac{\dot f^i-\dot f^j}{\lambda_i-\lambda_j}+ \frac{\dot f^i}{\lambda_j}+ \frac{\dot f^j}{\lambda_i}\geq 0 \]
	for each $i\neq j$ and each $\lambda\in \Gamma_+$ whose components are pairwise distinct. 
\end{lemma}

\bigskip

For $1$-homogeneous admissible underlying speed function, the local characterization of inverse-concavity is simplified. 

\begin{lemma}[{\cite[Lemmas~3 and 8]{andrews2013contracting}}]\label{inverse-concave_1-homogeneous}
	Let $f: \Gamma_+ \to \mathbb{R}$ be a function satisfying \textnormal{Conditions~\ref{conditions_f}(i)--(iii)}. Then $f$ is inverse-concave if and only if 
	\[\ddot f^{ij} +2 \frac{\dot f^i}{\lambda_i}\delta_{ij} \]
	is nonnegative definite for all $\lambda\in \Gamma_+$ and, in addition, 
	\[\frac{\dot f^i-\dot f^j}{\lambda_i-\lambda_j}+ \frac{\dot f^i}{\lambda_j}+ \frac{\dot f^j}{\lambda_i}\geq 0\]
	for each $i\neq j$ and each $\lambda\in \Gamma_+$ whose components are pairwise distinct. 
\end{lemma}

Note that the positivity follows from Conditions~\ref{conditions_f}(ii)(iii), together with the Euler relation $f(\lambda)=\sum_{i=1}^n \frac{\partial f}{\partial \lambda_i}(\lambda)\lambda_i$. 

\bigskip

\subsection{The inscribed and exscribed curvatures}\label{subsec:ballcurvatures}

For notational convenience, in this subsection, we shall not distinguish between a point of an embedded hypersurface and its image under the embedding.

\begin{definition}[Inscribed and exscribed curvatures, \cite{andrews2012noncollapsing, andrews2013non, sheng2009singularity}]
	Let $M^n = \partial \Omega^{n+1} \hookrightarrow \mathbb{R}^{n+1}$ be a smooth, properly embedded hypersurface bounding a precompact open set $\Omega^{n+1}\subset \mathbb{R}^{n+1}$, and let $\nu$ denote the outward unit normal vector field on $M^n$. 
	Define the extrinsic ball curvature $k: (M\times M)\setminus \{(x,x): x\in M\} \to \mathbb{R}$ by
	\[
	k(x,y)\coloneq \frac{2\left< x-y, \nu(x) \right>}{|x-y|^2}. 
	\]
	The inscribed curvature $\overline k: M \to \mathbb{R}$ is then given by
	\[
	\overline k(x) \coloneq \sup_{y\in M \setminus \{x\}} k(x,y), 
	\]
	and the exscribed curvature $\underline k: M \to \mathbb{R}$ by
	\[
	\underline k(x) \coloneq \inf_{y\in M \setminus \{x\}} k(x,y).
	\]
\end{definition}

\bigskip

To state the next lemma, which gives a geometric characterization of the ball curvatures, we adopt the following conventions. 
An \textit{extrinsic ball} $B$ will refer to either an open ball, an open half-space, or the complement of a closed ball in $\mathbb{R}^{n+1}$, equipped with its outward unit normal vector field $\nu_B$. In what follows, the curvature of $\partial B$ is understood to mean its constant principal curvature. In particular, it is taken to be positive if $B$ is a ball, and negative if $B$ is the complement of a ball.
We say that an extrinsic ball $B$ \textit{passes through} a point $y\in \mathbb{R}^{n+1}$ if $y\in \partial B$, and that $B$ \textit{touches} a hypersurface $M^n$ at $x\in M^n$ if $x\in \partial B$ and $\nu(x)=\nu_B(x)$.  

\begin{lemma}[\cite{andrews2016two, andrews2013non}]
	Let $M^n = \partial \Omega^{n+1} \hookrightarrow \mathbb{R}^{n+1}$ be a smooth, properly embedded hypersurface, equipped with its outward unit normal vector field $\nu$. Then $k(x,y)$ equals the curvature of the boundary of the unique extrinsic ball that touches the hypersurface $M$ at $x$ and passes through the point $y$. Moreover, $\overline k(x)$ equals the curvature of the boundary of the largest extrinsic ball contained in $\Omega$ that touches $M$ at $x$, while $\underline k(x)$ equals the curvature of the boundary of the smallest extrinsic ball containing $\Omega$ that touches $M$ at $x$.
\end{lemma}

Note that 
\[
\overline k(x)\geq \limsup_{y\to x} k(x,y)= \kappa_n(x)
\]
and 
\[
\underline k(x)\leq \liminf_{y\to x}k(x,y)= \kappa_1(x), 
\]
where, from this point on, we adopt the convention that the principal curvatures are ordered as $\kappa_1\leq \cdots\leq \kappa_n$, so that $\kappa_1$ and $\kappa_n$ denote the smallest and largest principal curvatures, respectively. 

\bigskip

We also recall the following lemma, which will be used in the proof of Corollary~\ref{cor}. 

\begin{lemma}[{\cite[Proposition~6.6]{langford2015motion}}]\label{lemma_langford}
	Let $X: M^n \to \mathbb{R}^{n+1}$ be an embedded convex hypersurface. Suppose that $\max_{M}\overline k$ is attained at some point $x_0\in M$. Then 
	\[\overline{k}(x_0) = \kappa_n(x_0).\]
\end{lemma}

\bigskip

\section{Auxiliary computations for the extrinsic ball curvature}\label{sec:aux}

In this section, we carry out preparatory computations concerning the evolution of the extrinsic ball curvature. The resulting formula~\eqref{evo_k} will be applied in the analysis of the ``interior case" in the subsequent proofs. 

Let us fix a point $(x_0,t_0) \in M^n \times (0,T)$ where $\underline k(x_0,t_0) < \kappa_1(x_0,t_0)$. Then there exists some $y_0 \in M^n \setminus\{x_0\}$ such that $\underline k(x_0,t_0) = k(x_0, y_0, t_0)$. We shall need to compute certain derivatives of $k$ at $(x_0,y_0,t_0)$. To this end, choose local normal coordinates $(x^1, \ldots, x^n): U_{x_0}\to \mathbb{R}^n$ and $(y^1, \ldots, y^n): U_{y_0} \to \mathbb{R}^n$ centered at $x_0$ and $y_0$, respectively, with respect to the metric at time~$t_0$. We can assume that $U_{x_0} \cap U_{y_0}=\emptyset$, so that $k$ is smooth in $U_{x_0}\times U_{y_0} \times (0,t_0]$.

To simplify notation, we define 
\[d(x,y,t) \coloneq |X(x,t)-X(y,t)|, \quad w(x,y,t) \coloneq \frac{X(x,t)-X(y,t)}{d(x,y,t)}, \]
and 
\[\partial_i^x(x,y,t) \coloneq \partial_{x^i}X (x,t), \quad \partial_i^y(x,y,t) \coloneq \partial_{y^i}X (y,t), \]
and use sub- and superscripts $x$ and $y$ to denote quantities relating to the first and second factors, respectively. 

With these notations in place, we compute
\[\partial_{x^i}k = -\frac{2}{d}(k \delta_i^p-{h^x}_i^p)\left< w, \partial_p^x \right> \]
and 
\[\partial_{y^i}k = -\frac{2}{d^2}\left< \partial_i^y, \nu_x-kdw \right>. \]
Since $\underline k$ is defined by taking the infimum over the second factor, we can use the vanishing of the $y$-derivatives at $y_0$ to determine the tangent plane to $M$ at $y_0$. 

\begin{lemma}
	Given an embedded hypersurface $M^n \subset \mathbb{R}^{n+1}$, suppose that $y_0 \in M^n \setminus \{x_0\}$ is an extremum of $k(x_0, \cdot)$. Then 
	\[\nu_y=\nu_x-kdw \]
	at $(x_0,y_0)$. 
\end{lemma}

Thus, the tangent plane at $y_0$ is the reflection of the tangent plane at $x_0$ across the plane orthogonal to $y_0-x_0$. In particular, we can choose the basis at $y_0$ to be the reflection of the basis at $x_0$. That is, 
\begin{equation}\label{relation_frames}
	\partial_i^y=\partial_i^x-2\left< \partial_i^x, w \right> w
\end{equation}
at $(x_0,y_0,t_0)$. 

We now differentiate once more to compute the second derivatives of $k$. Using the Codazzi equations, we derive at $(x_0,y_0,t_0)$: 
\[
\begin{aligned}
	\partial_{x^i} \partial_{x^j} k = & -k {h^x}_i^p {h^x}_{pj} + k^2 {h^x}_{ij} +\frac{4}{d^2}(k\delta_i^p-{h^x}_i^p)\left< w, \partial_j^x \right>\left< w, \partial_p^x \right> +\frac{4}{d^2}(k\delta_j^p-{h^x}_j^p)\left< w, \partial_i^x \right>\left< w, \partial_p^x \right> \\
	& -\frac{2}{d^2}(k\delta_{ij}-{h^x}_{ij}) + \frac{2}{d}\nabla^p {h^x}_{ij}\left< w,\partial_p^x \right>.
\end{aligned}
\]
Similarly, at $(x_0,y_0,t_0)$, we have
\[
\begin{aligned}
	\partial_{x^i} \partial_{y^j} k & = -\frac{4}{d^2}(k\delta_i^p-{h^x}_i^p)\left< w, \partial_p^x \right> \left< w, \partial_j^y \right> +\frac{2}{d^2}(k\delta_i^p-{h^x}_i^p)\left< \partial_p^x, \partial_j^y \right> \\
	& = \frac{2}{d^2}(k\delta_{ij}-{h^x}_{ij})
\end{aligned}
\]
and 
\[
\begin{aligned}
	\partial_{y^i} \partial_{y^j} k & = \frac{2}{d^2}\left< {h^y}_{ij}\nu_y, \nu_x-kdw \right>-\frac{2}{d^2}\left< \partial_i^y, k \partial_j^y \right> \\
	& = -\frac{2}{d^2}(k\delta_{ij}-{h^y}_{ij}).
\end{aligned}
\]
Combining these observations, we obtain at $(x_0,y_0,t_0)$: 
\[
\begin{aligned}
	& -\Psi'_x \dot F^{ij}_x (\partial_{x^i}+\Lambda_i^p \partial_{y^p})(\partial_{x^j}+\Lambda_j^q \partial_{y^q})k \\
	= & -\Psi'_x \dot F^{ij}_x (\partial_{x^i} \partial_{x^j}k+ 2\Lambda_i^p \partial_{x^j} \partial_{y^p}k + \Lambda_i^p \Lambda_j^q \partial_{y^p} \partial_{y^q}k) \\
	= & \ \Psi'_x |A^x|_F^2 k -\Psi'_x F_x k^2-\frac{8}{d^2}\Psi'_x \dot F^{ij}_x (k\delta_i^p-{h^x}_i^p)\left< w, \partial_j^x \right> \left< w, \partial_p^x \right>- \frac{2}{d} \left< w,\nabla \Psi_x \right> \\
	& + \frac{2}{d^2} \Psi'_x \dot F^{ij}_x \left[ (k\delta_{ij}-{h^x}_{ij})-2\Lambda_i^p(k\delta_{pj}-{h^x}_{pj})+ \Lambda_i^p\Lambda_j^q (k\delta_{pq}-{h^y}_{pq}) \right], 
\end{aligned}
\]
where $\Lambda$ is an arbitrary $n\times n$ matrix, to be determined in the next section. 

On the other hand, we compute the time derivative at $(x_0, y_0, t_0)$:
\[
\begin{aligned}
	\partial_t k & = \Psi_x k^2 -\frac{2}{d^2}\Psi_x +\frac{2}{d^2}\Psi_y\left< \nu_x-kdw, \nu_y \right> +\frac{2}{d}\left< w, \nabla \Psi_x \right> \\
	& = \Psi_x k^2 -\frac{2}{d^2}\Psi_x +\frac{2}{d^2}\Psi_y +\frac{2}{d}\left< w, \nabla \Psi_x \right>. 
\end{aligned}
\]

Thus, we conclude that 
\begin{equation}\label{evo_k}
	\begin{aligned}
		& \ \partial_t k -\Psi'_x \dot F^{ij}_x (\partial_{x^i}+\Lambda_i^p \partial_{y^p})(\partial_{x^j}+\Lambda_j^q \partial_{y^q})k \\
		= & \ \Psi'_x |A^x|_F^2 k -(\Psi'_x F_x-\Psi_x) k^2 -\frac{8}{d^2}\Psi'_x \dot F^{ij}_x (k\delta_i^p-{h^x}_i^p)\left< w, \partial_j^x \right> \left< w, \partial_p^x \right> \\
		& + \frac{2}{d^2} \left\{ \Psi_y-\Psi_x+ \Psi'_x \dot F^{ij}_x \left[ (k\delta_{ij}-{h^x}_{ij})-2\Lambda_i^p(k\delta_{pj}-{h^x}_{pj})+ \Lambda_i^p\Lambda_j^q (k\delta_{pq}-{h^y}_{pq}) \right] \right\}
	\end{aligned}
\end{equation}
at $(x_0,y_0,t_0)$. 

Observe that the first term on the right-hand side corresponds precisely to the reaction term, and the third term is positive at an off-diagonal $y$-minimum of $k$, where $k=\underline k<\kappa_1$. We will demonstrate that, under suitable conditions on the underlying speed function $f$ and the modulating function $\Psi$, the remaining terms can also be appropriately controlled. 

\bigskip

\section{Evolution of the exscribed curvature under the flow}\label{sec:evo_underlinek}

In this section, we establish the evolution equation for the exscribed curvature in the viscosity sense. This serves as an essential intermediate step toward the proofs of Theorems~\ref{thm_convex} and \ref{thm_inverse-concave}.

\begin{proposition}\label{prop_supersolution}
	Suppose that $f: \Gamma \to \mathbb{R}$ is an admissible underlying speed function and $\Psi$ satisfies \textnormal{Conditions~\ref{conditions_psi}(i)}. Let $X: M^n \times [0,T) \to \mathbb{R}^{n+1}$ be an embedded solution of \textnormal{\eqref{CF}}. Suppose also that either 
	\begin{itemize}
		\item[\rm{(i)}] $f$ is convex, and $\Psi$ satisfies \textnormal{Conditions~\ref{conditions_psi}(iiia)}, or
		\item[\rm{(ii)}] $\Gamma=\Gamma_+$ and $f$ is inverse-concave, and $\Psi$ satisfies \textnormal{Conditions~\ref{conditions_psi}(iv)}. 
	\end{itemize}
	Then the exscribed curvature $\underline k$ satisfies 
	\begin{equation*}
		(\partial_t-\Psi' \Delta_F)\underline k \geq \Psi' |A|_F^2 \underline k - (\Psi' F-\Psi)\underline k^2 
		\tag{$\ast$} \label{star}
	\end{equation*}
	in the viscosity sense.
\end{proposition}

\begin{proof}
	To keep the presentation streamlined, auxiliary lemmas will be introduced as needed and their verification postponed until the end of the argument. 
	
	Given a point $(x_0,t_0)\in M^n \times (0,T)$, suppose that $\phi \in C^\infty (U_{x_0}\times (t_0-\varepsilon ,t_0])$ is a lower support function for $\underline k$ at $(x_0,t_0)$; that is, 
	\[\phi \leq \underline k \ \ \text{in} \ \ U_{x_0}\times (t_0-\varepsilon ,t_0] \ \ \text{with} \ \ \phi(x_0,t_0)=\underline k(x_0,t_0). \]
	Then we need to prove that the differential inequality 
	\[(\partial_t-\Psi' \Delta_F)\phi \geq \Psi' |A|_F^2 \phi - (\Psi' F-\Psi)\phi^2\]
	holds at $(x_0,t_0)$. 
	
	By the definition of $\underline k$, there are two situations to consider: the ``interior case'', in which $\underline k(x_0,t_0)=k(x_0,y_0,t_0)$ for some $y_0\in M^n \setminus \{x_0\}$, and the ``boundary case'', in which $\underline k(x_0,t_0)=h_{(x_0,t_0)}(y_0,y_0)$ for some $y_0 \in S_{x_0}M$. 
	
	\bigskip

	We begin by considering the case where $f$ is convex. 
	
	\begin{enumerate}[nosep]
		\item[(i)] In this setting, we recall that $\Psi$ satisfies Conditions~\ref{conditions_psi}(iiia). 
	\end{enumerate}
	
	\medskip
	
	\noindent \textbf{The interior case}
	
	\noindent Observe that $\phi(x_0,t_0)=\underline k(x_0,t_0)=k(x_0,y_0,t_0)<\kappa_1(x_0,t_0)$ for some $y_0\in M \setminus \{x_0\}$ and $\phi(x,t)\leq \underline k(x,t)\leq k(x,y,t)$ for all $(x,t)\in U_{x_0}\times (t_0-\varepsilon, t_0]$ and all $y\in M \setminus \{x\}$. This implies that $\partial_t \phi(x_0,t_0)\geq \partial_t k(x_0,y_0,t_0)$, the gradient of $\phi-k$ on $M\times M$ vanishes at $(x_0,y_0,t_0)$, and the Hessian of $\phi-k$ on $M\times M$ is nonpositive definite at $(x_0,y_0,t_0)$. Thus, applying \eqref{evo_k}, we obtain at $(x_0,y_0,t_0)$: 
	\[
	\begin{aligned}
		(\partial_t-\Psi' \Delta_F)\phi \geq & \ \partial_t k- \Psi'_x \dot F^{ij}_x (\partial_{x^i}+\Lambda_i^p \partial_{y^p})(\partial_{x^j}+\Lambda_j^q \partial_{y^q})k \\
		= & \ \Psi'_x |A^x|_F^2 k -(\Psi'_x F_x-\Psi_x) k^2 -\frac{8}{d^2}\Psi'_x \dot F^{ij}_x (k\delta_i^p-{h^x}_i^p)\left< w, \partial_j^x \right> \left< w, \partial_p^x \right> \\
		& + \frac{2}{d^2} \left\{ \Psi_y-\Psi_x+ \Psi'_x \dot F^{ij}_x \left[ (k\delta_{ij}-{h^x}_{ij})-2\Lambda_i^p(k\delta_{pj}-{h^x}_{pj})+ \Lambda_i^p\Lambda_j^q (k\delta_{pq}-{h^y}_{pq}) \right] \right\}.
	\end{aligned}
	\]
	Taking $\Lambda$ to be the identity matrix yields 
	\[
	\begin{aligned}
		(\partial_t-\Psi' \Delta_F)\phi & \geq \Psi'_x |A^x|_F^2 \phi -(\Psi'_x F_x-\Psi_x) \phi^2 +\frac{2}{d^2}\left[ \Psi_y-\Psi_x+\Psi'_x \dot F^{ij}_x ({h^x}_{ij}-{h^y}_{ij}) \right] \\
		& \geq \Psi'_x |A^x|_F^2 \phi -(\Psi'_x F_x-\Psi_x) \phi^2,
	\end{aligned}
	\]
	where the last inequality follows from the convexity of $f$, which implies 
	\[\dot F^{ij}_x {h^y}_{ij}\leq F_y \]
	for any $x,y\in M$, together with Conditions~\ref{conditions_psi}(iiia) on $\Psi$, which ensures that 
	\[\Psi(b) - \Psi(a) +\Psi'(a)(a-b)\geq 0 \]
	for any $a,b\in (0,+\infty)$. This completes the proof of the interior case. 
	
	\bigskip

	\noindent \textbf{The boundary case}
	
	\noindent Here we have $\phi(x_0,t_0)=\underline k(x_0,t_0)=h_{(x_0,t_0)}(y_0,y_0)=\kappa_1(x_0,t_0)$ for some $y_0\in S_{x_0}M$. Let $\{e_i^0\}_{i=1}^n$ be a principal orthonormal frame at $(x_0,t_0)$ such that $y_0=e_1^0$. To proceed, we extend $\{e_i^0\}_{i=1}^n$ to a neighborhood of $(x_0,t_0)$ as follows: we first extend it spatially by parallel transport along radial geodesics emanating from $x_0$, with respect to the metric $g(t_0)$; then, we extend it in time by solving the ODE 
	\[\frac{d}{dt}e_i= \Psi \mathcal{W}(e_i) \]
	for each $i=1, \ldots, n$, where $\mathcal{W}$ is the Weingarten map. Then the resulting local frame field $\{e_i(t)\}_{i=1}^n$ remains orthonormal under the metric $g(t)$: 
	\[
	\begin{aligned}
		\frac{d}{dt}g_{ij} & = \left( \frac{d}{dt}g \right)_{ij}+ g\left(\frac{d}{dt}e_i, e_j\right)+ g\left( e_i, \frac{d}{dt}e_j \right) \\
		& = -2\Psi h_{ij} + \Psi h_{ij} +\Psi h_{ij} \\
		& = 0.
	\end{aligned}
	\]
	
	Define $u\coloneq h(e_1,e_1)$ in the neighborhood of $(x_0,t_0)$. Then, by applying the evolution equation for the second fundamental form, we obtain the following lemma: 
	\begin{lemma}\label{evo_u}
		The function $u$ satisfies 
		\[(\partial_t- \Psi' \Delta_F)u = (\Psi' \ddot F^{pq,rs}+ \Psi'' \dot F^{pq} \dot F^{rs}) \nabla_1 h_{pq} \nabla_1 h_{rs} +\Psi' |A|_F^2 u - (\Psi' F - \Psi)u^2 \]
		at $(x_0,t_0)$.
	\end{lemma}
	
	We now proceed with the proof of the boundary case. Since $\phi(x,t)\leq \underline k(x,t) \leq \kappa_1(x,t) \leq u(x,t)$ for all $(x,t)\in U_{x_0}\times (t_0-\varepsilon,t_0]$ and $\phi(x_0,t_0)=u(x_0,t_0)$, we deduce that at $(x_0,t_0)$, 
	\[
	\begin{aligned}
		(\partial_t-\Psi'\Delta_F)\phi & \geq (\partial_t-\Psi'\Delta_F)u \\
		& = (\Psi' \ddot F^{pq,rs}+ \Psi'' \dot F^{pq} \dot F^{rs}) \nabla_1 h_{pq} \nabla_1 h_{rs} +\Psi' |A|_F^2 u - (\Psi' F - \Psi)u^2 \\
		& \geq \Psi' |A|_F^2 \phi - (\Psi' F - \Psi)\phi^2,
	\end{aligned}
	\]
	where we have used the convexity of $f$ and Conditions~\ref{conditions_psi}(iiia) for $\Psi$ in the final step. This establishes the desired conclusion. 
	
	\bigskip

	We now turn to the case where $f$ is inverse-concave, which requires a slightly more delicate analysis. 
	
	\begin{enumerate}[nosep]
		\item[(ii)] In this setting, we recall that $\Psi$ satisfies Conditions~\ref{conditions_psi}(iv). 
	\end{enumerate}
	
	\medskip
	
	\noindent \textbf{The interior case}
	
	\noindent The proof proceeds in a manner similar to that of the convex case, and it suffices to prove the following lemma. 
	
	\begin{lemma}\label{inverse-concave_interior}
		Let $f: \Gamma_+ \to \mathbb{R}$ be an inverse-concave admissible underlying speed function, and suppose that $\Psi$ satisfies \textnormal{Conditions~\ref{conditions_psi}(i)(iv)}. Then, for any positive constant $k$, any positive definite matrix $A$ with $\lambda_{\min}(A)>k$, and any positive definite matrix $B$ with $\lambda_{\min}(B)\geq k$, the following inequality holds: 
		\[
		\begin{aligned}
			0\leq \Psi(F(B))- \Psi(F(A))+ \Psi'(F(A)) \dot F^{ij}(A) \sup_{\Lambda} \big[ & (k\delta_{ij}-A_{ij})-2\Lambda_i^p(k\delta_{pj}-A_{pj}) \\
			& + \Lambda_i^p\Lambda_j^q (k\delta_{pq}-B_{pq}) \big].
		\end{aligned}
		\]
	\end{lemma}
	
	\bigskip

	\noindent \textbf{The boundary case} 
	
	\noindent Recall that $\phi(x_0,t_0)= \underline k(x_0,t_0)= h_{(x_0,t_0)}(y_0,y_0)= \kappa_1(x_0,t_0)$ for some $y_0\in S_{x_0}M$. Let $\{x^i\}_{i=1}^n$ be geodesic normal coordinates for $(M,g(t_0))$ centered at $x_0$, and denote by $y^j=dx^j$ the corresponding fiber coordinates. Moreover, we may assume that $\{\partial_{x^i}|_{x_0}\}_{i=1}^n$ is a principal orthonormal frame at $(x_0,t_0)$, with $y_0=\partial_{x^1}|_{x_0}$. Define the function $K: TM\times [0,T) \to \mathbb{R}$ by $K(x,y,t)\coloneq h_{(x,t)}(y,y)$. Note that the function $\Phi(x,y,t) \coloneq \phi(x,t) g_{(x,t)}(y,y)$ is a lower support for $K$ at $(x_0,y_0,t_0)$. Therefore, at $(x_0,y_0,t_0)$, we have 
	\[\partial_t(K-\Phi) -\Psi'_x \dot F^{ij}_x(\partial_{x^i} - \Lambda_i^p \partial_{y^p})(\partial_{x^j} - \Lambda_j^q \partial_{y^q})(K-\Phi) \leq 0, \]
	where $\Lambda$ is an arbitrary $n\times n$ matrix, to be determined later. 
	
	We now proceed to compute the derivatives of $K-\Phi$ explicitly. Writing locally as $K-\Phi= (h_{kl}-\phi \, g_{kl})y^k y^l$, we compute 
	\[\partial_{x^i}(K-\Phi)= (\partial_{x^i}h_{kl}-\phi\, \partial_{x^i}g_{kl}-\partial_{x^i}\phi \, g_{kl})y^k y^l, \quad \partial_{y^i}(K-\Phi)= 2(h_{kl}-\phi \, g_{kl})y^k\delta_i^l. \]
	Evaluating these at $(x_0,y_0,t_0)$, we obtain 
	\begin{equation}\label{1st-order}
		\partial_{x^i}h_{11}-\partial_{x^i}\phi=0 \ \ \text{and} \ \ h_{1i}-\phi \, \delta_{1i}=0, \ \ \text{for} \ i=1, \ldots, n,
	\end{equation}
	since the first derivatives of $K-\Phi$ must vanish at $(x_0,y_0,t_0)$. We next compute the second derivatives: 
	\[
	\begin{aligned}
		\partial_{x^i}\partial_{x^j}(K-\Phi) & = (\partial_{x^i}\partial_{x^j}h_{kl}-\phi\, \partial_{x^i}\partial_{x^j}g_{kl}-\partial_{x^i}\phi \, \partial_{x^j}g_{kl}-\partial_{x^j}\phi \, \partial_{x^i}g_{kl}-\partial_{x^i}\partial_{x^j}\phi \, g_{kl})y^k y^l, \\
		\partial_{x^i}\partial_{y^j}(K-\Phi) & = 2(\partial_{x^i}h_{kl}-\phi\, \partial_{x^i}g_{kl}-\partial_{x^i}\phi \, g_{kl})y^k \delta_j^l, \\
		\partial_{y^i}\partial_{y^j}(K-\Phi) & = 2(h_{kl}-\phi \, g_{kl})\delta_i^k\delta_j^l.
	\end{aligned}
	\]
	Evaluating again at $(x_0,y_0,t_0)$, we find
	\[
	\begin{aligned}
		\partial_{x^i}\partial_{x^j}(K-\Phi) & = \partial_{x^i}\partial_{x^j}h_{11}-\phi\, \partial_{x^i}\partial_{x^j}g_{11}-\partial_{x^i}\partial_{x^j}\phi, \\
		\partial_{x^i}\partial_{y^j}(K-\Phi) & = 2(\partial_{x^i}h_{1j}-\partial_{x^i}\phi \, \delta_{1j}), \\
		\partial_{y^i}\partial_{y^j}(K-\Phi) & = 2(h_{ij}-\phi \, \delta_{ij}). 
	\end{aligned}
	\]
	Moreover, at $(x_0,t_0)$, the coordinate derivatives relate to the covariant derivatives via 
	\[
	\begin{aligned}
		\nabla_i \nabla_j h_{11} & = \partial_{x^i}\partial_{x^j}h_{11}-2\, \partial_{x^i}\Gamma_{1j}^k \, h_{1k}, \\
		0= \nabla_i \nabla_j g_{11} & = \partial_{x^i}\partial_{x^j}g_{11}-2\, \partial_{x^i}\Gamma_{1j}^k \, \delta_{1k}.
	\end{aligned}
	\]
	In view of the first-order condition~\eqref{1st-order}, it follows that 
	\[\partial_{x^i}\partial_{x^j}h_{11}-\phi\, \partial_{x^i}\partial_{x^j}g_{11}= \nabla_i \nabla_j h_{11}. \]
	Combining the above observations, we obtain at $(x_0,y_0,t_0)$: 
	\[
	\begin{aligned}
		& -\Psi'_x \dot F^{ij}_x(\partial_{x^i} - \Lambda_i^p \partial_{y^p})(\partial_{x^j} - \Lambda_j^q \partial_{y^q})(K-\Phi) \\
		= & - \Psi'_x \dot F^{ij}_x \left[ \partial_{x^i}\partial_{x^j}(K-\Phi)- 2\Lambda_i^p \partial_{x^j} \partial_{y^p}(K-\Phi) + \Lambda_i^p \Lambda_j^q \partial_{y^p} \partial_{y^q}(K-\Phi) \right] \\
		= & \ \Psi' \Delta_F \phi - \Psi' \Delta_F h_{11} + 2\Psi'\dot F^{ij} \Lambda_i^p \left[ 2(\nabla_j h_{1p}-\nabla_j h_{11}\delta_{1p}) - \Lambda_j^q(h_{pq}-h_{11} \delta_{pq}) \right]. 
	\end{aligned}
	\]
	We now turn to the time derivative. Using 
	\[\partial_t(K-\Phi)= (\partial_t h_{kl} -\phi\, \partial_t g_{kl} -\partial_t \phi \, g_{kl})y^k y^l, \]
	we find at $(x_0,y_0,t_0)$: 
	\[\partial_t(K-\Phi)= \partial_t h_{11}-h_{11}\partial_t g_{11}- \partial_t \phi. \]
	Substituting the evolution equations for the second fundamental form and the metric into the above expression, we obtain at $(x_0,y_0,t_0)$, 
	\[
	\begin{aligned}
		\partial_t(K-\Phi) = & -\partial_t \phi +\Psi'|A|_F^2 \phi -(\Psi' F -\Psi)\phi^2 \\
		& +\Psi' \Delta_F h_{11}+ (\Psi' \ddot F^{pq,rs}+ \Psi'' \dot F^{pq} \dot F^{rs}) \nabla_1 h_{pq} \nabla_1 h_{rs}.
	\end{aligned}  
	\]
	Thus, we conclude that 
	\[
	\begin{aligned}
		0\geq & \ \partial_t(K-\Phi) -\Psi'_x \dot F^{ij}_x(\partial_{x^i} - \Lambda_i^p \partial_{y^p})(\partial_{x^j} - \Lambda_j^q \partial_{y^q})(K-\Phi) \\
		= & -\partial_t \phi+ \Psi' \Delta_F \phi+ \Psi'|A|_F^2 \phi -(\Psi' F -\Psi)\phi^2+ (\Psi' \ddot F^{pq,rs}+ \Psi'' \dot F^{pq} \dot F^{rs}) \nabla_1 h_{pq} \nabla_1 h_{rs} \\
		& +  2\Psi'\dot F^{ij} \Lambda_i^p \left[ 2(\nabla_j h_{1p}-\nabla_j h_{11}\delta_{1p}) - \Lambda_j^q(h_{pq}-h_{11} \delta_{pq}) \right]
	\end{aligned}
	\]
	at $(x_0,y_0,t_0)$. Observe that the term in the square brackets vanishes when $p=1$. To arrive at the desired conclusion, we are left to prove that at $(x_0,t_0)$, 
	\[
	\begin{aligned}
		& (\Psi' \ddot F^{pq,rs}+ \Psi'' \dot F^{pq} \dot F^{rs}) \nabla_1 h_{pq} \nabla_1 h_{rs} \\
		& + 2\Psi'\dot F^{ij} \sup_{\Lambda} \left[ 2\Lambda_i^p(\nabla_j h_{1p}-\nabla_j h_{11}\delta_{1p}) - \Lambda_i^p \Lambda_j^q(h_{pq} -h_{11}\delta_{pq}) \right] \geq 0. 
	\end{aligned}
	\]
	
	We begin by employing an observation due to Brendle~\cite[Proposition~8]{brendle2013embedded} (see also~\cite[Lemma~6.15]{langford2015motion}): 
	\begin{lemma}\label{lemma_brendle}
		Let $X: M^n \to \mathbb{R}^{n+1}$ be an embedded hypersurface and suppose that $\underline k=\kappa_1$ at some point $x_0\in M^n$. Then $(\nabla_{e_1}h)(e_1,e_1)$ vanishes at $x_0$, where $e_1$ denotes the principal direction corresponding to $\kappa_1$.
	\end{lemma}
	
	Applying Lemma~\ref{lemma_brendle} and the following lemma, we complete the proof of the boundary case. 
	\begin{lemma}\label{inverse-concave_boundary}
		Let $f: \Gamma_+ \to \mathbb{R}$ be an inverse-concave admissible underlying speed function, and suppose that $\Psi$ satisfies \textnormal{Conditions~\ref{conditions_psi}(i)(iv)}. Then, for any diagonal matrix $A$ with positive eigenvalues $\lambda_1 \leq \cdots \leq \lambda_n$, and any symmetric matrix $B$ with $B_{11}=0$, the following inequality holds:
		\[
		\begin{aligned}
			& \ddot F^{pq,rs}(A)B_{pq}B_{rs}+ \frac{\Psi''(F(A))}{\Psi'(F(A))}\dot F^{pq}(A)\dot F^{rs}(A)B_{pq}B_{rs} \\
			& + 2\dot F^{ij}(A)\sup_{\Lambda_i^1=0}\left[ 2\Lambda_i^p B_{jp}- \Lambda_i^p \Lambda_j^q (A_{pq}-\lambda_1 \delta_{pq}) \right] \geq 0.
		\end{aligned} 
		\]
	\end{lemma}
	
	This concludes the proof of Proposition~\ref{prop_supersolution}.
\end{proof}

\bigskip

We now complete the argument by verifying the auxiliary Lemmas~\ref{evo_u}, \ref{inverse-concave_interior}, and \ref{inverse-concave_boundary} introduced along the way.

\begin{proof}[Proof of Lemma~\ref{evo_u}]
	A straightforward computation yields at $(x_0,t_0)$:
	\[
	\begin{aligned}
		\partial_t u & = (\partial_t h)_{11}+ 2h(\partial_t e_1, e_1) \\
		& = \Psi' \Delta_F h_{11}+ (\Psi' \ddot F^{pq,rs}+ \Psi'' \dot F^{pq} \dot F^{rs}) \nabla_1 h_{pq} \nabla_1 h_{rs} + \Psi' |A|_F^2 u - (\Psi' F - \Psi)u^2.
	\end{aligned} 
	\]
	It suffices to prove that 
	\[\Delta_F h_{11}=\Delta_F u\]
	at $(x_0,t_0)$. Noting that $\Gamma_{ij}^k=0$ at $(x_0,t_0)$ for all $i,j,k=1, \ldots, n$, we derive at $(x_0,t_0)$:
	\[
	\begin{aligned}
		\nabla_i \nabla_j h_{11} & = e_i((\nabla_{e_j}h)(e_1,e_1)) \\
		& = e_i(e_j(h(e_1,e_1)))-2e_i(h(\nabla_{e_j}e_1,e_1)) \\
		& = \nabla_i \nabla_j u - 2h(\nabla_{e_i}\nabla_{e_j}e_1,e_1).
	\end{aligned}
	\]
	Meanwhile, at $(x_0,t_0)$, we have
	\[\left< \nabla_{e_i}\nabla_{e_j}e_1,e_1 \right> = e_i \left< \nabla_{e_j}e_1, e_1 \right>=0, \]
	which in turn implies that $h(\nabla_{e_i}\nabla_{e_j}e_1,e_1)=0$ at $(x_0,t_0)$. Therefore, it follows that $\nabla_i \nabla_j h_{11}=\nabla_i \nabla_j u$ at $(x_0,t_0)$, and hence $\Delta_F h_{11}=\Delta_F u$. This completes the proof.
\end{proof}

\bigskip

\begin{proof}[Proof of Lemma~\ref{inverse-concave_interior}]
	By a standard approximation argument, we may assume $\lambda_{\min}(B)>k$, so that $B-kI$ is positive definite. Since the expression in the square brackets is quadratic in $\Lambda$, the supremum is attained for the choice $\Lambda=(A-kI)(B-kI)^{-1}$, where $I$ denotes the identity matrix. Consequently, it suffices to show that 
	\[0\leq \Psi(F(B))-\Psi(F(A))-\Psi'(F(A))\dot F^{ij}(A)\left[(A-kI)-(A-kI)(B-kI)^{-1}(A-kI) \right]_{ij}. \]
	
	To this end, define for $z\in [0,k]$ the function 
	\[q(z)\coloneq \Psi(F(B))-\Psi(F(A))-\Psi'(F(A))\dot F^{ij}(A)\left[(A-zI)-(A-zI)(B-zI)^{-1}(A-zI) \right]_{ij}. \]
	Differentiating gives
	\[q'(z)=\Psi'(F(A))\dot F^{ij}(A)\left( I-Y^{-1}X -XY^{-1}+XY^{-2}X \right)_{ij}, \]
	where we have set $X\coloneq A-zI$ and $Y\coloneq B-zI$. Now observe that 
	\[I-Y^{-1}X -XY^{-1}+XY^{-2}X= X(X^{-1}-Y^{-1})^2 X\geq 0, \]
	and therefore $q'(z)\geq 0$ for all $z\in [0,k]$. It follows that $q$ is non-decreasing. In particular, evaluating at the endpoints yields
	\begin{align}\label{q(0)}
		q(k) = \ & \Psi(F(B))-\Psi(F(A))-\Psi'(F(A))\dot F^{ij}(A)\left[(A-kI)-(A-kI)(B-kI)^{-1}(A-kI) \right]_{ij} \notag \\
		\geq \ & \Psi(F(B))-\Psi(F(A))-\Psi'(F(A))\dot F^{ij}(A)(A-A B^{-1} A)_{ij} =q(0) \notag \\
		= \ & \Psi(F(B))-\Psi(F(A))-\Psi'(F(A))F(A)+\Psi'(F(A))\dot F^{ij}(A)(AB^{-1}A)_{ij}.
	\end{align}
	At this stage, we invoke the inverse-concavity of $F$, which implies that 
	\[F(B)^{-1}=F_{*}(B^{-1})\leq \dot F_{*}^{ij}(A^{-1})(B^{-1})_{ij}=F(A)^{-2} \dot F^{ij}(A)(AB^{-1}A)_{ij}. \]
	Substituting this inequality into \eqref{q(0)}, we arrive at
	\[
	\begin{aligned}
		& \Psi(F(B))-\Psi(F(A))-\Psi'(F(A))F(A)+\Psi'(F(A))\dot F^{ij}(A)(AB^{-1}A)_{ij} \\
		\geq \ & \Psi(F(B))-\Psi(F(A))-\Psi'(F(A))F(A)+\Psi'(F(A))F(A)^2 F(B)^{-1}.
	\end{aligned}
	\]
	Thus, it remains only to verify that
	\[\Psi(b)-\Psi(a)-\Psi'(a)a^2\left(\frac{1}{a}-\frac{1}{b}\right)\geq 0, \]
	for any $a,b\in (0,+\infty)$. Indeed, by Conditions~\ref{conditions_psi}(iv), the function $s \mapsto \Psi'(s)s^2$ is non-decreasing on $(0,+\infty)$. Hence, if $a<b$, then for every $s\geq a$ we have $\Psi'(s)s^2\geq \Psi'(a)a^2$, or equivalently $\Psi'(s)\geq \Psi'(a)a^2/s^2$. Integrating this inequality over the interval $[a,b]$ yields
	\[\Psi(b)-\Psi(a)\geq \Psi'(a)a^2\left(\frac{1}{a}-\frac{1}{b}\right). \] 
	The case $a>b$ follows by the same argument. In conclusion, the inequality holds for all $a,b \in (0,+\infty)$, which establishes the claim and completes the proof. 
\end{proof}

\bigskip

\begin{proof}[Proof of Lemma~\ref{inverse-concave_boundary}]
	It suffices to prove the claim in the case where all $\lambda_i$ are distinct, since any diagonal matrix can be approximated by such matrices. 
	
	Set
	\[
	\begin{aligned}
		Q\coloneq & \ \ddot F^{pq,rs}(A)B_{pq}B_{rs}+ \frac{\Psi''(F(A))}{\Psi'(F(A))}\dot F^{pq}(A)\dot F^{rs}(A)B_{pq}B_{rs} \\
		& + 2\dot F^{ij}(A)\sup_{\Lambda_i^1=0}\left[ 2\Lambda_i^p B_{jp}- \Lambda_i^p \Lambda_j^q (A_{pq}-\lambda_1 \delta_{pq}) \right].
	\end{aligned} 
	\]
	Note that the supremum over $\Lambda$ is attained when $\Lambda_i^p= (\lambda_p-\lambda_1)^{-1}B_{ip}$ for $p>1$. With this choice, we obtain 
	\[Q= \ddot F^{pq,rs}(A)B_{pq}B_{rs}+ \frac{\Psi''(F(A))}{\Psi'(F(A))}\dot F^{pq}(A)\dot F^{rs}(A)B_{pq}B_{rs}+ 2\dot F^{ij}(A) R^{pq}B_{ip}B_{jq}, \]
	where $R^{pq}\coloneq (\lambda_p-\lambda_1)^{-1}\delta^{pq}$ for $p,q\neq 1$, and zero otherwise. To reach the desired result, we only need to verify that 
	\[\left(\ddot F^{pq,rs}(A)+ \frac{\Psi''(F(A))}{\Psi'(F(A))}\dot F^{pq}(A)\dot F^{rs}(A)+ 2\dot F^{pr}(A) R^{qs} \right) B_{pq}B_{rs}\geq 0 \]
	for any symmetric matrix $B$ with $B_{11}=0$. Applying~\cite[Theorem~5.1]{andrews2007pinching}, we deduce that 
	\[
	\begin{aligned}
		& \left(\ddot F^{pq,rs}(A)+ \frac{\Psi''(F(A))}{\Psi'(F(A))}\dot F^{pq}(A)\dot F^{rs}(A)+ 2\dot F^{pr}(A) R^{qs} \right) B_{pq}B_{rs} \\
		= & \, \sum_{p,q}\ddot f^{pq}(\lambda)B_{pp}B_{qq}+ \sum_{p\neq q} \frac{\dot f^p(\lambda)-\dot f^q(\lambda)}{\lambda_p-\lambda_q}B_{pq}^2 \\
		& + \frac{\Psi''(f(\lambda))}{\Psi'(f(\lambda))}\sum_{p,q}\dot f^p(\lambda) \dot f^q(\lambda) B_{pp}B_{qq}+ 2\sum_{p=1}^n \sum_{q=2}^n\frac{\dot f^p(\lambda)}{\lambda_q-\lambda_1}B_{pq}^2 \\
		= & \, \sum_{p,q}\ddot f^{pq}(\lambda)B_{pp}B_{qq}+ \frac{\Psi''(f(\lambda))}{\Psi'(f(\lambda))}\sum_{p,q}\dot f^p(\lambda) \dot f^q(\lambda) B_{pp}B_{qq}+ 2\sum_{p,q=2}^n \frac{\dot f^p(\lambda) \delta_{pq}}{\lambda_p-\lambda_1}B_{pp}B_{qq} \\
		& + \sum_{p\neq q} \frac{\dot f^p(\lambda)-\dot f^q(\lambda)}{\lambda_p-\lambda_q}B_{pq}^2+ 2\sum_{p=2}^n \frac{\dot f^1(\lambda)}{\lambda_p-\lambda_1}B_{1p}^2+ 2\sum_{\substack{p,q>1 \\p\neq q}}\frac{\dot f^p(\lambda)}{\lambda_q-\lambda_1}B_{pq}^2. 
	\end{aligned}
	\]
	We first estimate 
	\[
	\begin{aligned}
		& \, \sum_{p,q}\ddot f^{pq}(\lambda)B_{pp}B_{qq}+ \frac{\Psi''(f(\lambda))}{\Psi'(f(\lambda))}\sum_{p,q}\dot f^p(\lambda) \dot f^q(\lambda) B_{pp}B_{qq}+ 2\sum_{p,q=2}^n \frac{\dot f^p(\lambda) \delta_{pq}}{\lambda_p-\lambda_1}B_{pp}B_{qq} \\
		\geq & \, \sum_{p,q}\ddot f^{pq}(\lambda)B_{pp}B_{qq}-\frac{2}{f(\lambda)}\sum_{p,q}\dot f^p(\lambda) \dot f^q(\lambda) B_{pp}B_{qq}+ 2\sum_{p,q=2}^n \frac{\dot f^p(\lambda) \delta_{pq}}{\lambda_p}B_{pp}B_{qq} \\
		= & \, \sum_{p,q} \left(\ddot f^{pq}(\lambda)- \frac{2}{f(\lambda)}\dot f^p(\lambda) \dot f^q(\lambda)+ 2\frac{\dot f^p(\lambda)}{\lambda_p}\delta_{pq} \right)B_{pp}B_{qq}\geq 0, 
	\end{aligned}
	\]
	where we have used Conditions~\ref{conditions_psi}(iv) for $\Psi$ in the first step, and the final inequality follows from the inverse-concavity of $f$ (Lemma~\ref{inverse-concave_origin}). The remaining terms can be estimated as follows: 
	\[
	\begin{aligned}
		& \, \sum_{p\neq q} \frac{\dot f^p(\lambda)-\dot f^q(\lambda)}{\lambda_p-\lambda_q}B_{pq}^2+ 2\sum_{p=2}^n \frac{\dot f^1(\lambda)}{\lambda_p-\lambda_1}B_{1p}^2+ 2\sum_{\substack{p,q>1 \\p\neq q}}\frac{\dot f^p(\lambda)}{\lambda_q-\lambda_1}B_{pq}^2 \\
		= & \, \sum_{\substack{p,q>1 \\p\neq q}}\left(\frac{\dot f^p(\lambda)-\dot f^q(\lambda)}{\lambda_p-\lambda_q}+2\frac{\dot f^p(\lambda)}{\lambda_q-\lambda_1} \right)B_{pq}^2+ 2\sum_{p=2}^n \left( \frac{\dot f^p(\lambda)-\dot f^1(\lambda)}{\lambda_p-\lambda_1} + \frac{\dot f^1(\lambda)}{\lambda_p-\lambda_1} \right)B_{1p}^2 \\
		\geq & \, \sum_{\substack{p,q>1 \\p\neq q}}\left(\frac{\dot f^p(\lambda)-\dot f^q(\lambda)}{\lambda_p-\lambda_q}+\frac{\dot f^p(\lambda)}{\lambda_q}+\frac{\dot f^q(\lambda)}{\lambda_p} \right)B_{pq}^2+ 2\sum_{p=2}^n \frac{\dot f^p(\lambda)}{\lambda_p-\lambda_1} B_{1p}^2.
	\end{aligned}
	\]
	The first term is nonnegative due to the inverse-concavity of $f$ (Lemma~\ref{inverse-concave_origin} or Lemma~\ref{inverse-concave_1-homogeneous}), while the second term is clearly nonnegative. This completes the proof.
\end{proof}

\bigskip

\section{Proof of Theorem~\ref{thm_convex}}\label{sec:proof_convex}

In this section, we give the proof of Theorem~\ref{thm_convex}. 
\begin{proof}
	Consider the function $Z(x,t)\coloneq \underline k(x,t) +\beta F(x,t)$, and it suffices to prove that $Z$ remains nonnegative under the flow. Given any $\tau \in [0,T)$ and $\varepsilon>0$, we consider the function $Z_{\varepsilon}\coloneq Z+\varepsilon e^{(C+1)t}$, where $C\coloneq \max_{M \times [0,\tau]} \Psi'|A|_F^2 <+\infty$. Observe that $Z_{\varepsilon}$ is positive at the initial time. We claim that $Z_{\varepsilon}$ remains positive under the flow up to time $\tau$. Since $\tau$ and $\varepsilon$ are arbitrary, this will suffice to prove the claim. Suppose, to the contrary, that there exists some $(x_0,t_0)\in M \times (0,\tau]$ such that $Z_{\varepsilon}(x_0,t_0)=0$. We may assume that $t_0$ is the first such time, and this means that the function $\phi(x,t)\coloneq -\beta F(x,t)-\varepsilon e^{(C+1)t}$ is a lower support for $\underline k$ at $(x_0,t_0)$. Since $\underline k$ is a viscosity supersolution of \eqref{star}, we obtain that at $(x_0,t_0)$, $\phi$ satisfies 
	\[
	\begin{aligned}
		0 \geq & -(\partial_t-\Psi' \Delta_F)\phi + \Psi' |A|_F^2 \phi - (\Psi' F-\Psi) \phi^2 \\
		= & \ \varepsilon (C+1)e^{(C+1)t_0} + \beta(|A|_F^2 \Psi + \Psi''|\nabla F|_F^2) \\
		& - \Psi'|A|_F^2 (\beta F+\varepsilon e^{(C+1)t_0})-(\Psi' F-\Psi)(\beta F+\varepsilon e^{(C+1)t_0})^2 \\
		\geq & \ \varepsilon (C+1)e^{(C+1)t_0}-\varepsilon \Psi' |A|_F^2 e^{(C+1)t_0}+ \beta(\Psi-\Psi'F)|A|_F^2 \\
		& + \beta \Psi'' |\nabla F|_F^2 +(\Psi-\Psi' F)(\beta F+\varepsilon e^{(C+1)t_0})^2 \\
		\geq & \ \varepsilon (C+1)e^{(C+1)t_0}-\varepsilon C e^{(C+1)t_0} \\
		= & \ \varepsilon e^{(C+1)t_0}>0,
	\end{aligned}
	\]
	which is a contradiction. This completes the proof. 
\end{proof}

\bigskip

\section{Proof of Theorem~\ref{thm_inverse-concave} and Corollary~\ref{cor}}\label{sec:proof_inverse-concave}

This section contains the proofs of Theorem~\ref{thm_inverse-concave} and Corollary~\ref{cor}.

\medskip

We begin by proving Theorem~\ref{thm_inverse-concave}.
\begin{proof}[Proof of Theorem~\ref{thm_inverse-concave}]
	We first observe that 
	\[0<\underline k F \leq \kappa_1 \sum_{i=1}^n \frac{\partial f}{\partial \lambda_i}(\kappa)\kappa_i \leq \sum_{i=1}^n \frac{\partial f}{\partial \lambda_i}(\kappa)\kappa_i^2 =|A|_F^2. \]
	Combining it with Conditions~\ref{conditions_psi}(iib) yields
	\[-(\Psi'F-\Psi)\underline k^2 \geq -\frac{\Psi'F-\Psi}{F}|A|_F^2 \underline k.\]
	It follows that the exscribed curvature $\underline k$ satisfies 
	\begin{equation*}
		(\partial_t-\Psi' \Delta_F)\underline k \geq \frac{\Psi}{F} |A|_F^2 \underline k
		\tag{$\ast'$} \label{star'}
	\end{equation*}
	in the viscosity sense. Since the underlying speed function $F$ satisfies \eqref{evo_F}, the statement of Theorem~\ref{thm_inverse-concave} follows from a simple comparison argument for viscosity solutions: 
	
	Define $u(t)\coloneq \min_{x\in M}\frac{\underline k(x,t)}{F(x,t)}$ for each time $t$. We claim that $u(t)$ is non-decreasing in $t$. It suffices to prove that 
	\[\underline k(\cdot, t)- (u(t_0)-\varepsilon e^{t-t_0})F(\cdot, t) \geq 0 \]
	for any $t_0\in [0,T), \ t\in [t_0,T)$ and $\varepsilon>0$. Taking $\varepsilon \to 0$ then gives $\underline k(\cdot, t)-u(t_0)F(\cdot, t)\geq 0$, and hence $u(t)\geq u(t_0)$ for $t\geq t_0$. 
	
	Fix $t_0\in [0, T)$ and $\varepsilon>0$. At time $t_0$, we have $\underline k(\cdot, t_0)-(u(t_0)-\varepsilon) F(\cdot, t_0)\geq \varepsilon F(\cdot, t_0)>0$. Suppose, to the contrary, that there is some $(x_1,t_1)\in M \times (t_0,T)$ such that $\underline k(x_1,t_1)-(u(t_0)-\varepsilon e^{t_1-t_0}) F(x_1,t_1)=0$. We may assume that $t_1$ is the first such time, which means that the function $\phi(x,t) \coloneq (u(t_0)-\varepsilon e^{t-t_0}) F(x,t)$ is a lower support for $\underline k$ at $(x_1,t_1)$. Since $\underline k$ is a viscosity supersolution of \eqref{star'}, we obtain that at $(x_1,t_1)$, $\phi$ satisfies 
	\[
	\begin{aligned}
		0 & \geq -(\partial_t-\Psi' \Delta_F)\phi+\frac{\Psi}{F}|A|_F^2\phi \\
		& = \varepsilon e^{t_1-t_0} F- (u(t_0)-\varepsilon e^{t_1-t_0})(|A|_F^2 \Psi + \Psi''|\nabla F|_F^2) +\frac{\Psi}{F}|A|_F^2 (u(t_0)-\varepsilon e^{t_1-t_0})F \\
		& = \varepsilon e^{t_1-t_0} F - (u(t_0)-\varepsilon e^{t_1-t_0})\Psi''|\nabla F|_F^2 \\
		& \geq \varepsilon e^{t_1-t_0}F>0,
	\end{aligned}
	\]
	where the relation $u(t_0)-\varepsilon e^{t_1-t_0}=\frac{\underline k(x_1,t_1)}{F(x_1,t_1)}>0$ is applied at the penultimate step. This contradiction implies that $\underline k(\cdot, t)- (u(t_0)-\varepsilon e^{t-t_0})F(\cdot, t)$ does not vanish for any $t\in [t_0,T)$. This completes the proof. 
\end{proof}

\bigskip

We then deduce Corollary~\ref{cor}. 
\begin{proof}[Proof of Corollary~\ref{cor}]
	By Theorem~\ref{thm_inverse-concave}, and recalling that $\underline k(x,t)$ denotes the curvature of the smallest enclosing ball of $M_t$ touching the hypersurface at $X(x,t)$, we obtain  
	\[\frac{1}{\rho_+(t)}\geq \max_{x\in M}\underline k(x,t)\geq c_0 \max_{x\in M} F(x,t), \]
	where $c_0=\min_{M \times \{0\}}\frac{\underline k}{F}>0$. 
	Moreover, by the monotonicity of $F$, there holds 
	\[F\geq \gamma \kappa_1, \]
	with $\gamma \coloneq f(1, \cdots, 1)$. 
	
	On the other hand, the argument in the ``boundary case" also yields 
	\[\frac{\kappa_1}{F}\geq \min_{M\times \{0\}} \frac{\kappa_1}{F}\eqcolon k_0>0. \]
	Consequently, along the flow~\eqref{CF}, the $n$-tuple of principal curvatures remains in the cone 
	\[\Gamma_0 \coloneq \bigg\{(\lambda_1, \dots, \lambda_n)\in \Gamma_+ : \frac{\min\{\lambda_1, \dots, \lambda_n\}}{f(\lambda_1, \dots, \lambda_n)}\geq k_0\bigg\}. \]
	Under the additional assumption that
	\[\frac{\min\{\lambda_1, \dots, \lambda_n\}}{f(\lambda_1, \dots, \lambda_n)} \text{ approaches zero on the boundary of } \Gamma_+,\]
	it follows immediately that $\Gamma_0 \Subset \Gamma_+$. We show that, when $n=2$, this conclusion holds without imposing the above condition. Using the $1$-homogeneity of $f$, we may write
	\[f(\lambda_1,\lambda_2) = \lambda_2 f(\lambda_1/\lambda_2, 1). \]
	Assuming $\lambda_1\leq \lambda_2$ and setting $s=\lambda_1/\lambda_2$, we consider 
	\[v(\lambda_1,\lambda_2)\coloneq \frac{\min\{\lambda_1, \lambda_2\}}{f(\lambda_1,\lambda_2)} = \frac{s}{f(s,1)}.\]
	Thus $v$ depends only on the ratio $s$. Differentiating, we obtain  
	\[v'(s) = \frac{f(s,1)-s\frac{\partial f}{\partial \lambda_1}(s,1)}{f^2(s,1)}. \]
	By the Euler relation,
	\[f(\lambda_1,\lambda_2) = \lambda_1 \frac{\partial f}{\partial \lambda_1}(\lambda_1,\lambda_2)+ \lambda_2 \frac{\partial f}{\partial \lambda_2}(\lambda_1,\lambda_2), \]
	and evaluating at $(s,1)$ gives
	\[v'(s)=\frac{\frac{\partial f}{\partial \lambda_2}(s,1)}{f^2(s,1)}>0. \]
	Hence $v$ is monotone increasing in $s$. Recall that $k_0=\min_{M\times\{0\}}\frac{\kappa_1}{F}$, which is equivalently the minimum of $v$ over a compact set $K\subset \Gamma_+$. Choose a cone $C_0$ such that $K \subset C_0\Subset \Gamma_+$. The monotonicity of $v$ in $s$ then implies $\Gamma_0 \subset C_0$, and therefore $\Gamma_0 \Subset \Gamma_+$.
	As a consequence, 
	\[\kappa_1\geq \delta \kappa_n, \]
	where 
	\[\delta\coloneq \inf_{\Gamma_0}\frac{\min\{\lambda_1, \dots, \lambda_n\}}{\max\{\lambda_1, \dots, \lambda_n\}}>0.\] 
	
	Finally, let $x_0\in M$ be a point at which $\kappa_n(\cdot,t)$ attains its maximum. By Lemma~\ref{lemma_langford}, the function $\overline k(\cdot,t)$ also attains its maximum at $x_0$, and $\overline k(x_0,t)= \kappa_n(x_0,t)$. This implies that there exists a ball of radius $1/\kappa_n(x_0,t)$ which is enclosed by $M_t$ and touches the hypersurface at $X(x_0,t)$. In particular, 
	\[\frac{1}{\rho_-(t)}\leq \kappa_n(x_0,t)= \max_{x\in M}\kappa_n(x,t). \]
	Combining the above estimates yields 
	\[\rho_+(t)\leq \frac{1}{\gamma \delta c_0}\rho_-(t), \]
	which completes the proof. 
\end{proof}

\medskip

\begin{remark*}
	For $n\geq 3$, one cannot expect $\Gamma_0 \Subset \Gamma_+$ without additional assumptions. We illustrate this with an explicit example in the case $n=3$. 
	
	Consider the harmonic mean 
	\[f(\lambda_1, \lambda_2, \lambda_3) = \frac{3}{\frac{1}{\lambda_1}+\frac{1}{\lambda_2}+\frac{1}{\lambda_3}}.\]
	Assuming $\lambda_1 \leq \lambda_2 \leq \lambda_3$, we compute
	\[v(\lambda_1,\lambda_2,\lambda_3) \coloneq \frac{\min\{\lambda_1,\lambda_2,\lambda_3\}}{f(\lambda_1,\lambda_2,\lambda_3)} = \frac{1}{3}\bigg(1+\frac{\lambda_1}{\lambda_2}+\frac{\lambda_1}{\lambda_3}\bigg). \]
	
	Let $K\subset \Gamma_+$ denote the compact set of principal curvature triples of the initial hypersurface~$M_0$, and suppose that $(1,10,10)\in K$. Then 
	\[k_0 = \min_{K}v\leq v(1,10,10) = \frac{2}{5}.\]
	
	Now consider the sequence on the unit sphere
	\[\xi^{(m)}= \frac{1}{\sqrt{1^2+5^2+m^2}}(1,5,m). \]
	For $m$ sufficiently large, 
	\[v(\xi^{(m)}) = v(1,5,m) = \frac{2}{5}+\frac{1}{3m} > \frac{2}{5} \geq k_0,\]
	so that $\xi^{(m)} \in \Gamma_0$. On the other hand, 
	\[\lim_{m\to \infty} \xi^{(m)} = (0,0,1) \in \partial\Gamma_+, \]
	and hence $\Gamma_0 \not\Subset \Gamma_+$. 
\end{remark*}

\bigskip

\noindent \textbf{Acknowledgements.}
The second author would like to thank Mat Langford for useful discussions. 
Research of the first author was partially supported by the National Key R\&D Program of China (No. 2022YFA1005500) and the National Natural Science Foundation of China under Grant Nos. 12031017 and 12571063.


\end{document}